\documentclass{article}

\usepackage{amsthm}
\usepackage{amssymb}
\usepackage{amsmath}

\newtheorem{theorem}{Theorem}[section]

\newtheorem{lemma}[theorem]{Lemma}
\theoremstyle{remark}

\theoremstyle{definition}

\theoremstyle{definition}

\begin{document}

\title{Necessary Optimality Conditions
for a Dead Oil Isotherm Optimal Control Problem\footnote{Work
supported by the \emph{Portuguese Foundation for Science and
Technology} (FCT) through the \emph{Centre for Research in
Optimization and Control} (CEOC) of the University of Aveiro,
cofinanced by the European Community fund FEDER/POCTI. The first
author was also supported by the postdoc fellowship
SFRH/BPD/20934/2004. Research Report CM06/I-42.}}

\author{Moulay Rchid Sidi Ammi\footnote{Postdoc Researcher,
Department of Mathematics, University of Aveiro, 3810-193 Aveiro,
Portugal. Email: \texttt{sidiammi@mat.ua.pt}}
\and
Delfim F. M. Torres\footnote{Associate Professor,
Department of Mathematics, University of Aveiro,
3810-193 Aveiro, Portugal. Email: \texttt{delfim@mat.ua.pt}}}

\date{}

\maketitle


\begin{abstract}
We study a system of nonlinear partial differential equations
resulting from the traditional modelling of oil engineering within
the framework of the mechanics of a continuous medium. Recent
results on the problem provide existence, uniqueness and
regularity of the optimal solution. Here we obtain the first
necessary optimality conditions.
\end{abstract}


\smallskip


\smallskip


\smallskip


\medskip


\section{Introduction}

We are interested in the optimal control of the dead oil
isotherm problem:
\begin{equation}
\label{P}
\begin{cases}
\partial_t u - \Delta \varphi(u) = div\left(g(u) \nabla p\right)
& \text{ in } Q_T = \Omega \times (0,T) \, , \\
\partial_t p - div\left(d(u) \nabla p\right) = f
& \text{ in } Q_T = \Omega \times (0,T) \, , \\
\left.u\right|_{\partial \Omega} = 0 \, , \quad \left.u\right|_{t=0} = u_0 \, , \\
\left.p\right|_{\partial \Omega} = 0 \, , \quad \left.p\right|_{t=0} = p_0 \, ,
\end{cases}
\end{equation}
where $\Omega$ is an open bounded domain in $\mathbb{R}^2$ with a
sufficiently smooth boundary. Equations \eqref{P} serve as a model
for an incompressible biphasic flow in a porous medium, with
applications to the industry of exploitation of hydrocarbons. To
understand the optimal control problem we will consider here, some
words about the recovery of hydrocarbons are in order. At the time
of the first run of a layer, the flow of the crude oil towards the
surface is due to the energy stored in the gases under pressure
in the natural hydraulic system. To mitigate the consecutive
decline of production and the decomposition of the site, water
injections are carried out, well before the normal exhaustion of
the layer. The water is injected through wells with high pressure,
by pumps specially drilled to this end. The pumps allow the
displacement of the crude oil towards the wells of production. The
wells must be judiciously distributed, which gives rise to a
difficult problem of optimal control: how to choose the best
installation sites of the production wells? This is precisely the
question we address in this work: our main goal is to present a
method to carry out the optimal control of \eqref{P} with respect
to all the important parameters that intervene in the process.
More precisely, we seek necessary conditions for the admissible
parameters $u$, $p$ and $f$ to minimize the functional
\begin{equation}
\label{eq:cf} J(u,p,f) = \frac{1}{2} \left\|u -
U\right\|_{2,Q_T}^2 + \frac{1}{2} \left\|p - P\right\|_{2,Q_T}^2 +
\frac{\beta_1}{2} \left\|f\right\|_{2 q_0,Q_T}^{2 q_0} +
\frac{\beta_2}{2} \left\|\partial_t f\right\|_{2,Q_T}^{2}
\end{equation}
where $q_0 > 1$ and $\beta_1 > 0$ and $\beta_2 > 0$ are two
coefficients of penalization. The first two terms in \eqref{eq:cf}
make possible to minimize the difference between the reduced
saturation of oil $u$, the global pressure $p$ and the given data
$U$ and $P$.

Existence and uniqueness to the system \eqref{P}, for the case when
the term $\partial_t p$ is missing but for more general boundary
conditions, is established in \cite{mad}. In \cite{romaniaDG-DS} we
obtain conditions which provide existence and regularity of the
optimal solutions to the problem of minimizing \eqref{eq:cf} subject
to \eqref{P}. Here we are interested to obtain necessary optimality
conditions which permit to find the solutions predicted by the
results in \cite{romaniaDG-DS}. This is, to the best of our
knowledge, an important open question.

Several techniques for deriving optimality conditions are available
in the literature of optimal control systems governed by partial
differential equations \cite{Lions,Lions71,Boris2006}. We obtain
the optimality conditions by making use of a Lagrangian approach
recently used with success by O.~Bodart, A.~V.~Boureau and
R.~Touzani for an optimal control problem of the induction heating
\cite{Bodart}, and by H.-C.~Lee and T.~Shilkin for the thermistor
problem \cite{LeeShilkin}.


\section{Preliminaries}

Our main objective is to obtain necessary conditions for a triple
$\left(\bar{u},\bar{p},\bar{f}\right)$ to minimize \eqref{eq:cf}
among all the functions $\left(u,p,f\right)$ verifying \eqref{P}.
The intended necessary optimality conditions are proved in
\S\ref{sec:MR} under adequate hypotheses on the data of the
problem, which provide regularity of the optimal solution.

\subsection{Notation and Functional Spaces}

In the sequel we assume that $\varphi$, $g$ and $d$ are real
valued $C^1$-functions satisfying:

\begin{description}

\item[(H1)] $0 < c_1 \le d(r)$, $\varphi(r) \le c_2$;
$|d'(r)|,\, |\varphi'(r)|,\, |\varphi''(r)| \leq c_{3}\,  \quad
\forall r \in \mathbb{R}$.

\item[(H2)] $u_0$, $p_0$ $\in C^2\left(\bar{\Omega}\right)$, $U$,
$P$ $\in L^2(Q_T)$, where $u_0,\, p_0 : \Omega \rightarrow
\mathbb{R}$, $U,\, P : Q_T \rightarrow \mathbb{R}$, and
$\left.u_0\right|_{\partial \Omega} = \left.p_0\right|_{\partial
\Omega} = 0$.

\end{description}

We consider the following spaces:
\begin{equation*}
W_p^{1,0}(Q_T) := L^p\left(0,T,W_p^1(\Omega)\right)
= \left\{ u \in L^p(Q_T), \, \nabla u \in L^p(Q_T) \right\} \, ,
\end{equation*}
endowed with the norm
$\left\|u\right\|_{W_p^{1,0}(Q_T)}
= \left\|u\right\|_{p,Q_T} + \left\|\nabla u\right\|_{p,Q_T}$;
\begin{equation*}
W_p^{2,1}(Q_T) :=
\left\{ u \in W_p^{1,0}(Q_T), \, \nabla^2 u \, , \partial_t u \in L^p(Q_T) \right\} \, ,
\end{equation*}
with the norm $\left\|u\right\|_{W_p^{2,1}(Q_T)}
= \left\|u\right\|_{W_p^{1,0}(Q_T)} + \left\|\nabla^2 u\right\|_{p,Q_T}
+ \left\|\partial_t u\right\|_{p,Q_T}$;
\begin{gather*}
V := \left\{ u \in W_2^{1,0}(Q_T), \, \partial_t u
\in L^2\left(0,T,W_2^{-1}(\Omega)\right) \right\} \, ; \\
W := \left\{ u \in W_{2q}^{2,1}(Q_T), \, \left.u\right|_{\delta_T
= \partial \Omega \times (0,T)} = 0\right\} \, ; \\
\Upsilon := \left\{ f \in L^{2 q}(Q_T) , \, \partial_t f \in L^2(Q_T) \right\} \, ; \\
H := L^{2 q}(Q_T) \times
\stackrel{\circ}{W}_{2 q}^{2 - \frac{1}{q}}(\Omega) \, .
\end{gather*}


\subsection{Coercive Estimate}

The following lemma provides a coercive estimate to linear
parabolic systems that is useful for our purposes.

\begin{lemma}[\cite{sol}]
\label{lemma2.4} Let $\Omega$ be a bounded domain with a
$C^2$-boundary, and assume that
\begin{equation*}
A_{i j k l} \in C\left(\bar{Q}_T\right) \, , \quad
b_{i j k} \in L^r(Q_T) \, , \quad
c_{i j} \in L^{\frac{r}{2}}(Q_T)
\end{equation*}
with $r > n + 2$ and $A$ satisfying the strong ellipticity
condition:
\begin{gather*}
\exists \gamma_0 > 0 : A_{i j k l}(u) B_{i j} B_{k l} \ge \gamma_0 |B|^2
\quad \forall B \in M^{n \times n} \, , \\
A_{i j k l} A_{k l i j} = A_{j i k l} = A_{i j l k} \, .
\end{gather*}
Then, for any $s \in (1,r)$, $s \ne \frac{3}{2}$, and for
arbitrary functions $f \in L^s\left(Q_T,\mathbb{R}^n\right)$ and
$u_0 \in
\stackrel{\circ}{W}_s^{2-\frac{2}{s}}\left(\Omega,\mathbb{R}^N\right)$
there exists a unique solution $u \in W_s^{2,1}(Q_T,\mathbb{R}^N)$
of the problem
\begin{gather*}
\partial_t u_i - A_{i j k l}(x) u_{k j l}
+ b_{i j k}(x) u_{j k} + c_{i j}(x) u_j = f_i(x) \, ,\\
\left. u \right|_{t=0} = u_0 \, ,
\quad \left. u \right|_{\partial \Omega} = 0 \, .
\end{gather*}
Moreover, the estimate
\begin{equation*}
\left\|u\right\|_{W_s^{2,1}(Q_T)} \le c
\left(\left\|f\right\|_{s,Q_T} +
\left\|u_0\right\|_{W_s^{2-\frac{2}{s}}(\Omega)}\right)
\end{equation*}
holds for some constant $c$ depending only on $n$, $\Omega$, $T$,
$\gamma_0$ and the norms of the coefficients.
\end{lemma}


\subsection{Existence of Optimal Solution}

The following existence theorem is proved in \cite{romaniaDG-DS}
using a technical lemma found in \cite{blp}, Young's inequality
and Aubin's Lemma, together with the theorem of Lebesgue and some
compacity arguments of J. L. Lions \cite{Lions}. The conclusion
follows from the fact that $J$ is lower semicontinuous with
respect to the weak convergence.

\begin{theorem}[\cite{romaniaDG-DS}]
\label{theorem3.1} Under the hypotheses (H1) and (H2) there exists a
$q > 1$, depending on the data of the problem, such that the
problem of minimizing \eqref{eq:cf} subject to \eqref{P} has an
optimal solution $\left(\bar{u},\bar{p},\bar{f}\right)$ satisfying
$\bar{u} \in W_{q}^{2,1}(Q_T)$,
$\bar{p} \in C\left([0,T];L^2(\Omega)\right) \cap W_{2 q}^{1, 0}(Q_T)$,
$\partial_t \bar{p} \in L^2\left(0, T, W_2^{-1}(\Omega)\right)$,
$\bar{f} \in L^{2 q_0}(Q_T)$,
$\partial_t \bar{f} \in L^2(Q_T)$.
\end{theorem}


\subsection{Regularity of Solutions}

Regularity of the solutions given by Theorem~\ref{theorem3.1} is
also proved in \cite{romaniaDG-DS}. Theorem~\ref{theorem4.1} is
obtained using Young's and Holder's inequalities, Gronwall Lemma,
De Giorgi–-Nash–-Ladyzhenskaya–-Uraltseva theorem,
an estimate from \cite{ks}, and some technical lemmas found in \cite{lsu}.

\begin{theorem}[\cite{romaniaDG-DS}]
\label{theorem4.1} Let $\left(\bar{u},\bar{p},\bar{f}\right)$ be
an optimal solution to the problem of minimizing \eqref{eq:cf}
subject to \eqref{P}. Suppose that (H1) and (H2) are satisfied.
Then, there exist $\alpha > 0$ such that the following regularity
conditions hold:
$\bar{u} \in C^{\alpha,\frac{\alpha}{2}}\left(\bar{Q}_T\right)$,
$\bar{u}, \, \bar{p} \in W_{4}^{1,0}(Q_T)$,
$\bar{u}, \, \bar{p} \in W_{2}^{2,1}(Q_T)$,
$\partial_t \bar{u}, \, \partial_t \bar{p}
\in L^\infty\left(0,T; L^2(\Omega)\right) \cap W_2^{1,0}(Q_T)$,
$\bar{u} \in C^{\frac{1}{4}}(Q_T)$,
$\bar{u} \in W_{2 q_0}^{2,1}(Q_T)$,
$\bar{p} \in W_{2 q_0}^{2,1}(Q_T)$.
\end{theorem}


\section{Main Results}
\label{sec:MR}

We define the following nonlinear operator
corresponding to \eqref{P}:
\begin{gather*}
F : W \times W \times \Upsilon \longrightarrow H \times H \\
\left(u,p,f\right) \longrightarrow F(u,p,f) = 0
\end{gather*}
where
\begin{equation*}
F(u,p,f) =
\left(%
\begin{array}{cc}
  \partial_t u - \Delta \varphi(u) - div(g(u) \nabla p), & \gamma_0 u - u_0 \\
  \partial_t p - div\left(d(u) \nabla p\right) - f, & \gamma_0 p - p_0 \\
\end{array}%
\right) \, ,
\end{equation*}
$\gamma_0$ being the trace operator $\gamma_0 u =
\left.u\right|_{t=0}$. Owing to the estimate
\begin{equation*}
\left\|v\right\|_{W_{\frac{4 q}{2 - q}}^{1,0}}
\le c \left\|v\right\|_{W_{2 q}^{2,1}(Q_T)} \, ,
\quad \forall v \in W_{2 q}^{2,1}(Q_T) \, , \quad 1 < q < 2
\end{equation*}
(see \cite{lsu}) and hypothesis (H1), we have
\begin{equation*}
\varphi''(u) \left|\nabla u\right|^2 , \,
g'(u) \nabla u \nabla p , \,
d(u) \nabla u \nabla p \in L^{\frac{2 q}{2 - q}}(Q_T) \subset L^{2 q}(Q_T) \, .
\end{equation*}
Thus, it follows that $F$ is well defined.


\subsection{G\^{a}teaux differentiability}

\begin{theorem}
\label{thm5.1} In addition to the hypotheses (H1) and (H2), let us
suppose that
\begin{description}
\item[(H3)]
 $\left|\varphi'''\right| \le c$.
\end{description}
 Then, the operator $F$
is G\^{a}teaux differentiable and its derivative is given by
\begin{multline*}
\delta F(u,p,f)(e,w,h) = \frac{d}{ds} F\left(u + s e, p + s w, f + s
h\right)\left.\right|_{s=0}
= \left(\delta F_1, \delta F_2\right) \\
= \left(%
\begin{array}{cc}
  \partial_t e - div\left(\varphi'(u)\nabla e\right)
  - div\left(\varphi''(u) e \nabla u\right)
  - div\left(g(u) \nabla w\right) - div\left(g'(u) e \nabla p\right), & \gamma_0 e \\
  \partial_t w - div\left(d(u) \nabla w\right) - div\left(d'(u) e \nabla p\right) - h, & \gamma_0 w \\
\end{array}%
\right) \, ,
\end{multline*}
for all $(e,w,h) \in W \times W \times \Upsilon. $\\
 Furthermore, for any optimal solution
$\left(\bar{u},\bar{p},\bar{f}\right)$ of the problem of minimizing
\eqref{eq:cf} among all the functions $\left(u,p,f\right)$
satisfying \eqref{P}, the image of $\delta
F\left(\bar{u},\bar{p},\bar{f}\right)$ is equal to $H \times H$.
\end{theorem}

To prove Theorem~\ref{thm5.1} we make use of the following lemma.

\begin{lemma}
\label{lemma5.2} The operator $\delta F(u,p,f) : W \times W \times
\Upsilon \longrightarrow H \times H$ is  linear and bounded.
\end{lemma}

\begin{proof} of Lemma~\ref{lemma5.2}:
We have for all $(e,w,h) \in W \times W \times \Upsilon $
\begin{equation*}
\begin{split}
&\delta_{p}F_{2}(u, p, f)(e, w, h)= \partial_{t}w - div\left( d(u)
\nabla w\right)- div\left( d'(u)e \nabla p \right)-h\\
&= \partial_{t}w- d(u)\triangle w -d'(u) \nabla u. \nabla w   -
d'(u)e \triangle p -d'(u) \nabla e. \nabla u- d'(u)e \nabla u.
\nabla p -h \, ,
\end{split}
\end{equation*}
with $\delta_{p} F$ the G\^ateaux derivative of $F$ with respect to
$p$.
Then, using hypothesis $(H1)$, we obtain that
\begin{multline}
\label{eq:4.12}
\|\delta_{p}F_{2}(u, p, f)(e, w, h) \|_{2q, Q_{T}}\leq
\|\partial_{t}w\|_{2q,
Q_{T}}+
 \| \nabla w \|_{2q, Q_{T}}+ c \| \triangle w \|_{2q, Q_{T}}\\
 + c  \|\nabla u. \nabla w\|_{2q, Q_{T}}
 +c \| e \triangle p \|_{2q, Q_{T}}+ c  \|\nabla e. \nabla u
\|_{2q, Q_{T}}
+ c  \| e \nabla u. \nabla p\|_{2q, Q_{T}}+ \|h\|_{2q, Q_{T}} \, .
\end{multline}
In what follows we consider the term
$\| e \nabla u. \nabla p\|_{2q, Q_{T}}$.
Similar arguments apply to the remaining terms of \eqref{eq:4.12}.
We have
\begin{equation*}
\begin{split}
\| e \nabla u. \nabla p\|_{2q, Q_{T}} &\leq \|e\|_{\infty, Q_{T}} \|
\nabla u. \nabla p\|_{2q, Q_{T}}\\
&\leq  \|e\|_{\infty, Q_{T}} \|\nabla u\|_{\frac{4q}{2-q}, Q_{T}} \|\nabla
p\|_{4, Q_{T}}\\
&\leq c \|u\|_{W} \|p\|_{W} \|e\|_{W}\, .
\end{split}
\end{equation*}
Then,
\begin{equation}\label{eq:4.13}
\|\delta_{p}F_{2}(u, p, f)(e, w, h) \|_{2q, Q_{T}}
\leq c \left( \|u\|_{W}, \|p\|_{W}, \|f\|_{\Upsilon} \right)
\left( \|e\|_{W} +\|w\|_{W}+ \|h\|_{\Upsilon} \right)\, .
\end{equation}
On the other hand,
\begin{multline*}
\delta_{u}F_{1}(u, p, f)(e, w, h)
= \partial_{t}e -
div\left( \varphi'(u) \nabla e\right)- div\left( \varphi''(u)e
\nabla u \right)\\
- div\left( g(u) \nabla w\right)- div\left( g'(u)e
\nabla p \right)\\
= \partial_{t}e- \varphi'(u)\triangle e
-\varphi''(u) \nabla u. \nabla e  \\
 -\varphi''(u)e \triangle u -\varphi''(u) \nabla e. \nabla u
- \varphi'''(u)e |\nabla u|^{2}\\
-g(u)\triangle w- g'(u) \nabla u. \nabla w\\
- g'(u) e \triangle p- g'(u)\nabla e. \nabla p -g''(u)e \nabla u.
\nabla p \, ,
\end{multline*}
where $\delta_{u} F$ is the G\^ateaux derivative of $F$ with respect
to $u$. The same argument as above give that
\begin{equation}\label{eq:4.14}
\|\delta_{u}F_{1}(u, p, f)(e, w, h) \|_{2q, Q_{T}}
\leq c \left( \|u\|_{W}, \|p\|_{W}, \|f\|_{\Upsilon} \right)
\left( \|e\|_{W}+ \|w\|_{W}+ \|h\|_{\Upsilon} \right)\, .
\end{equation}
Consequently, by \eqref{eq:4.13} and \eqref{eq:4.14} we can write
\begin{equation*}
\|\delta F(u, p, f)(e, w, h) \|_{H\times H\times \Upsilon}
\leq c \left( \|u\|_{W}, \|p\|_{W}, \|f\|_{\Upsilon} \right)
\left( \|e\|_{W}+ \|w\|_{W}+ \|h\|_{\Upsilon} \right)\, .
\end{equation*}
\end{proof}
\begin{proof} of Theorem~\ref{thm5.1}:
In order to show that the image of $\delta F(\overline{u},
\overline{p}, \overline{f})$ is equal to $H \times H$, we need to
prove that there exists a $(w, e, h) \in W \times W \times \Upsilon$
such that
\begin{equation}
\label{eq:4.15}
\begin{gathered}
 \partial_{t}e -
div\left( \varphi'(\overline{u}) \nabla e\right)-
div\left( \varphi''(\overline{u})e
\nabla \overline{u} \right)
- div\left( g(\overline{u}) \nabla w\right)- div\left( g'(\overline{u})e
\nabla \overline{p} \right)= \alpha \, ,  \\
\partial_{t}w - div\left( d(\overline{u})
\nabla w\right)- div\left( d'(\overline{u})e
 \nabla \overline{p} \right)-h=
\beta  \, , \\
\left.e\right|_{\partial \Omega} = 0 \,
 , \quad \left.e\right|_{t=0} = b\, , \\
\left.w\right|_{\partial \Omega} = 0 \, ,  \quad
\left.w\right|_{t=0} = a \, ,
\end{gathered}
\end{equation}
for any $(\alpha, a)$ and $(\beta, b) \in H$. Writing the system
\eqref{eq:4.15} for $h=0$ as
\begin{equation}
\label{eq:4.16}
\begin{gathered}
 \partial_{t}e -
\varphi'(\overline{u})\triangle e
-2 \varphi''(\overline{u})\nabla \overline{u}. \nabla e-
\varphi''(\overline{u})e \triangle \overline{u}-
\varphi'''(\overline{u})e |\nabla \overline{u}|^{2} \, , \\
-g(\overline{u}) \triangle w-
g'(\overline{u}) \nabla \overline{u}. \nabla w-
g'(\overline{u}) e \triangle \overline{p}
-g'(\overline{u})\nabla \overline{p}. \nabla e
- g''(\overline{u}) e \nabla \overline{u}. \nabla \overline{p}
= \alpha \, ,  \\
\partial_{t}w -  d(\overline{u}) \triangle w
-d'(\overline{u}) \nabla \overline{u}. \nabla w
- d'(\overline{u}) e \triangle \overline{p} -
d'(\overline{u}) \nabla \overline{u}. \nabla \overline{e}
-d'(\overline{u}) e \nabla \overline{u}. \nabla \overline{p}
=\beta  \, ,\\
\left.e\right|_{\partial \Omega} = 0 \,
 , \quad \left.e\right|_{t=0} = b\, , \\
\left.w\right|_{\partial \Omega} = 0 \, ,  \quad
\left.w\right|_{t=0} = a \, ,
\end{gathered}
\end{equation}
it follows from the regularity of the optimal solution
(Theorem~\ref{theorem4.1}) that
\begin{gather*}
\varphi''(\overline{u}) \triangle \overline{u}\, ,
 \varphi'''(\overline{u}) |\nabla \overline{u}|^{2}\, ,
g'(\overline{u})  \triangle \overline{p}\, ,
g''(\overline{u})  \nabla \overline{u}. \nabla \overline{p}\, ,
d'(\overline{u})  \triangle \overline{p}\, ,
d'(\overline{u})  \nabla \overline{u}. \nabla \overline{p}
 \in L^{2q_{0}(Q_{T})} \, , \\
 \varphi''(\overline{u})\nabla \overline{u}\, ,
  g'(\overline{u})\nabla \overline{u} \, ,
   g'(\overline{u})\nabla \overline{p}\, ,
   d'(\overline{u})\nabla \overline{u}
   \in L^{4q_{0}}(Q_{T}) \, .
\end{gather*}
By Lemma~\ref{lemma2.4} there exists a unique solution of
the system \eqref{eq:4.16}, hence there exists a $(e, w, 0)$
verifying \eqref{eq:4.15}. We conclude that the image of
$\delta F$ is equal to $H \times H$.
\end{proof}


\subsection{Necessary Optimality Condition}

We consider the cost functional $J: W\times  W \times \Upsilon
\rightarrow \mathbb{R}$ \eqref{eq:cf} and the Lagrangian $\mathcal{L}$
defined by
 $$
 \mathcal{L}\left(u, p, f, p_{1}, e_{1}, a, b\right)=
 J\left(u, p, f \right)+ \left\langle F(u, p, f),
 \left(\begin{array}{cc}  p_{1} & a \\
                  e_{1}, & b
\end{array}\right) \right\rangle\, ,
$$
where the bracket $\langle \cdot, \cdot \rangle $ denote the duality
between $H$ and $H'$.
\begin{theorem}\label{thm5.2}
Under hypotheses (H1)--(H3), if $\left(\overline{u}, \overline{p},
\overline{f}\right)$ is an optimal solution to the problem of
minimizing \eqref{eq:cf} subject to \eqref{P}, then there exist
functions $(\overline{e_{1}}, \overline{p_{1}}) \in
W_{2}^{2,1}(Q_{T}) \times W_{2}^{2,1}(Q_{T})$ satisfying the
following conditions:
\begin{equation*}
\begin{gathered}
\partial_{t}\overline{e_{1}} +
div\left( \varphi'(\overline{u}) \nabla e_{1}\right)
 -d'(\overline{u}) \nabla \overline{p}. \nabla \overline{p_{1}}
-\varphi''(\overline{u}) \nabla \overline{u}. \nabla
\overline{e_{1}} - g'(\overline{u})\nabla \overline{p}. \nabla
\overline{e_{1}}= \overline{u}-U\, ,  \\
\left.\overline{e_{1}}\right|_{\partial \Omega} = 0 \,
 , \quad \left.\overline{e_{1}}\right|_{t=T} = 0\, ,
\end{gathered}
\end{equation*}
\begin{equation*}
\begin{gathered}
\partial_{t}\overline{p_{1}} + div\left( d(\overline{u})
\nabla \overline{p_{1}}\right)+ div\left( g(\overline{u}) \nabla
\overline{e_{1}}\right)=\overline{p}-P \, , \\
\left.\overline{p_{1}}\right|_{\partial \Omega} = 0 \, , \quad
\left.\overline{p_{1}}\right|_{t=T} = 0 \, ,
\end{gathered}
\end{equation*}
\begin{equation}
\label{eq:4.19} - \beta_{2}\frac{\partial^{2}
\overline{f}}{\partial t^{2}}
 +2q_{0} \beta_{1}|\overline{f}|^{2q_{0}-2}\overline{f}
  = \overline{p_{1}} \, , \quad \,
\left.\frac{\partial \overline{f}}{\partial t}\right|_{t=0} =
\left.\frac{\partial \overline{f}}{\partial t}\right|_{t=T}=0.
\end{equation}
\end{theorem}
\begin{proof}
Let $ \left( \overline{u}, \overline{p}, \overline{f} \right)$ be an
optimal solution to the problem of
minimizing \eqref{eq:cf} subject to \eqref{P}. It is well known
(\textrm{cf. e.g.} \cite{fur}) that there
exist Lagrange multipliers $\left( (\overline{p_{1}}, \overline{a}),
(\overline{e_{1}}, \overline{b}) \right) \in H' \times H'$ verifying
$$ \delta_{(u, p, f)}\mathcal{L}
\left(\overline{u}, \overline{p}, \overline{f}, \overline{p_{1}},
\overline{e_{1}}, \overline{a}, \overline{b} \right)(e, w, h)= 0 \,
\quad \forall (e, w, h)\in W \times W \times \Upsilon, $$ with
$\delta_{(u, p, f)}\mathcal{L}$ the G\^{a}teaux derivative of
$\mathcal{L}$ with respect to $(u, p, f)$. We then obtain
\begin{gather*}
\int_{Q_{T}} \left( (\overline{u}-U)e +(\overline{p}-P)w+ 2q_{0}
\beta_{1}|\overline{f}|^{2q_{0}-2}\overline{f}h +
\beta_{2}\partial_{t} \overline{f} \partial_{t} h \right)\, dx dt\\
-\int_{Q_{T}} \left( \partial_{t}e - div\left(
\varphi'(\overline{u}) \nabla e \right)- div\left(
\varphi''(\overline{u})e \nabla \overline{u}\right)- div\left(
g(\overline{u}) \nabla w \right) - div\left( g'(\overline{u}) e
\nabla \overline{p}\right)
\right)\overline{e_{1}} \, dx dt \\
-\int_{Q_{T}} \left( \partial_{t}w - div\left( d(\overline{u})
\nabla w \right)- div\left( d'(\overline{u}) e\nabla
\overline{p}\right)-h
\right)\overline{p_{1}}  \, dx dt \\
-\langle \gamma_{0} e, \overline{a}\rangle + -\langle \gamma_{0} w,
\overline{b}\rangle =0 \, \quad \forall (e, w, h)\in W \times W
\times \Upsilon.
\end{gather*}
This last system is equivalent to the following one:
\begin{equation*}
\begin{gathered}
\int_{Q_{T}} \left( (\overline{u}-U)e-div\left( d'(\overline{u})
e\nabla \overline{p}\right)\overline{p_{1}}+\partial_{t}e \,
\overline{e_{1}} -
 div\left(
\varphi'(\overline{u}) \nabla e \right) \overline{e_{1}} \right. \\
 \left. -div\left( \varphi''(\overline{u})e \nabla
\overline{u}\right) \overline{e_{1}} -div\left( g'(\overline{u}) e
\nabla \overline{p}\right)\overline{e_{1}} \right) \,
 dx dt  \\
 +\int_{Q_{T}} \left( (\overline{p}-P)w +\partial_{t}w \,
\overline{p_{1}}- div\left( d(\overline{u}) \nabla w \right)
\overline{p_{1}} - div\left( g(\overline{u}) \nabla w
\right)\overline{e_{1}} \right) \, dx dt  \\
 +\int_{Q_{T}} \left(2q_{0}
\beta_{1}|\overline{f}|^{2q_{0}-2}\overline{f}h +
\beta_{2}\partial_{t} \overline{f} \partial_{t} h- \overline{p_{1}}
h \right)\, dx dt  \\
 + \langle \gamma_{0} e, \overline{a}\rangle +
\langle \gamma_{0} w, \overline{b}\rangle =0 \, \quad \forall (e, w,
h)\in W \times W \times \Upsilon.
\end{gathered}
\end{equation*}
In others words, we have
\begin{equation}
\label{eq:4.21}
\begin{gathered}
\int_{Q_{T}} \left( (\overline{u}-U)+d'(u)\nabla \overline{p}.
\nabla \overline{p_{1}}- \partial_{t}\overline{ e_{1}}-
 div\left( \varphi'(\overline{u}) \nabla \overline{e_{1}}
\right) + \varphi''(\overline{u}) \nabla \overline{u}. \nabla
\overline{e_{1}}+g'(u)\nabla \overline{p}. \nabla \overline{e_{1}}
\right)e \, dx dt  \\
+\int_{Q_{T}} \left( (\overline{p}-P) +\partial_{t}\overline{p_{1}}
- div\left( d(\overline{u}) \nabla \overline{p_{1}} \right) -
div\left( g(\overline{u}) \nabla \overline{e_{1}} \right) \right)w
\, dx dt  \\
+\int_{Q_{T}} \left(2q_{0}
\beta_{1}|\overline{f}|^{2q_{0}-2}\overline{f}h +
\beta_{2}\partial_{t} \overline{f} \partial_{t} h- \overline{p_{1}}
h \right)\, dx dt  \\
+ \langle \gamma_{0} e, \overline{a}\rangle + \langle \gamma_{0} w,
\overline{b}\rangle =0 \,
 \forall (e, w, h)\in W \times W \times \Upsilon.
\end{gathered}
\end{equation}
Consider now the system
\begin{equation}
\label{eq:4.22}
\begin{gathered}
\partial_{t}e_{1}
+div\left( \varphi'(\overline{u}) \nabla e_{1} \right)-
d'(\overline{u})\nabla \overline{p}. \nabla p_{1}
-\varphi''(\overline{u}) \nabla \overline{u}. \nabla
e_{1}-g'(\overline{u})\nabla \overline{p}. \nabla e_{1}=
\overline{u}-U  \, ,\\
\partial_{t}p_{1}+div\left( d(\overline{u}) \nabla p_{1} \right)
+div\left( g(\overline{u}) \nabla e_{1} \right)= \overline{p}-P, \\
 \left.e_{1}\right|_{\partial \Omega}
 =\left.p_{1}\right|_{\partial \Omega}  = 0 \, ,
\quad \left.e_{1}\right|_{t=T} = \left.p_{1}\right|_{t=T} = 0 \, .
\end{gathered}
\end{equation}
It follows by Lemma~\ref{lemma2.4} that \eqref{eq:4.22} has a
unique solution $(e_{1}, p_{1})\in W_{2}^{2, 1}(Q_{T}) \times
W_{2}^{2,1}(Q_{T})$. Since the problem of finding $(e, w)\in W \times W$ satisfying
\begin{equation}
\label{eq:4.23}
\begin{gathered}
\partial_{t}e
-div\left( \varphi'(\overline{u}) \nabla e \right)- div\left(
\varphi''(\overline{u}) e \nabla \overline{u} \right)- div\left(
g(\overline{u}) \nabla w \right)- div\left( g'(\overline{u})e \nabla
\overline{p} \right)= sign(e_{1}- \overline{e_{1}})\,  \\
\partial_{t}w -div\left(d(\overline{u}) \nabla w \right)-
div\left( d'(\overline{u})e \nabla \overline{p} \right)=
sign(p_{1}- \overline{p_{1}})\,  \\
\gamma_{0}e= \gamma_{0}w= 0
\end{gathered}
\end{equation}
is uniquely solvable on $W_{2q}^{2, 1}\times W_{2q}^{2, 1}$ by
Lemma~\ref{lemma2.4}, choosing $h=0$ in \eqref{eq:4.21}, multiplying
\eqref{eq:4.22} by $(e, w)$, integrating by parts, and making the
difference with \eqref{eq:4.21}, we obtain
\begin{equation}
\label{eq:4.24}
\begin{gathered}
\int_{Q_{T}} \left(\partial_{t}e -div\left( \varphi'(\overline{u})
\nabla e \right)- div\left( \varphi''(\overline{u}) e \nabla
\overline{u} \right)- div\left( g(\overline{u}) \nabla w \right)-
div\left( g'(\overline{u})e \nabla \overline{p} \right)
\right)(e_{1}-
\overline{e_{1}})\, dx dt  \\
 +\int_{Q_{T}}\left( \partial_{t}w -div\left(d(\overline{u}) \nabla w \right)-
div\left( d'(\overline{u})e \nabla \overline{p} \right)
\right)(p_{1}- \overline{p_{1}})\, dx dt  \\
+\langle \gamma_{0}e,  \gamma_{0}\overline{e_{1} } - \overline{a}
\rangle
 + \langle \gamma_{0}w,  \gamma_{0}\overline{p_{1} } - \overline{b}
\rangle =0 \, \quad \forall (e, w) \in W \times W.
\end{gathered}
\end{equation}
Choosing $(e, w)$ in \eqref{eq:4.24} as the solution of the system
\eqref{eq:4.23}, we have
 \begin{gather*}
\int_{Q_{T}} sign(e_{1}- \overline{e_{1}})(e_{1}-
\overline{e_{1}})\, dx dt
 + \int_{Q_{T}} sign(p_{1}- \overline{p_{1}})(p_{1}-
\overline{p_{1}})\, dx dt = 0 \, .
\end{gather*}
It follows that $e_{1}=\overline{ e_{1}}$ and  $p_{1}=\overline{
p_{1}}$. Coming back to \eqref{eq:4.24}, we obtain
$\gamma_{0}\overline{e_{1}}= \overline{a}$ and
$\gamma_{0}\overline{p_{1}}= \overline{b}$. On the other hand,
choosing $(e, w)=(0, 0)$ in \eqref{eq:4.21} it
follows \eqref{eq:4.19}, which conclude the proof of Theorem~\ref{thm5.2}.
\end{proof}



\end{document}